\documentclass[12pt,twoside]{amsart}
\usepackage{amssymb}
\usepackage{amscd}
\usepackage[abbrev,alphabetic]{amsrefs}
\usepackage{hyperref}
\usepackage{amsmath}
\usepackage{xypic}
\usepackage[all,cmtip]{xy}
\usepackage{tikz}
\usepackage{tikz-cd}
\usepackage{caption}

\usepackage{url}
\RequirePackage{booktabs, multirow}
\RequirePackage{pgf}
\newcommand{\git}{\mathbin{
  \mathchoice{/\mkern-4mu/}% \displaystyle
    {/\mkern-5mu/}% \textstyle
    {/\mkern-5mu/}% \scriptstyle
    {/\mkern-5mu/}}}% 

\title[KV vanishing fails for log del Pezzo surfaces in char. $3$]
{Kawamata-Viehweg vanishing fails for log del Pezzo surfaces in characteristic 3} 
\author{Fabio Bernasconi} 

\address{Department of Mathematics, University of Utah, Salt Lake City, UT 84112, USA} 
\email{fabio@math.utah.edu}

\newcommand{\Spec}[0]{{\operatorname{Spec}}}

\newcommand{\Supp}[0]{{\operatorname{Supp}}}
\newcommand{\Pic}[0]{{\operatorname{Pic}}}
\newcommand{\car}[0]{{\operatorname{char}}}

\newcommand{\Exc}[0]{{\operatorname{Exc}}}

\newcommand{\Cl}[0]{{\operatorname{Cl}}}
\newtheorem{thm}{Theorem}[section]
\newtheorem{lem}[thm]{Lemma}
\newtheorem{cor}[thm]{Corollary}
\newtheorem{proposition}[thm]{Proposition}

\theoremstyle{definition}

\newtheorem{definition}[thm]{Definition}

\newtheorem{remark}[thm]{Remark}

\makeatletter
  
  \@addtoreset{equation}{section}
  \makeatother

\newcommand{\MO}{\mathcal{O}}

\newcommand{\Q}{\mathbb{Q}}
\newcommand{\Z}{\mathbb{Z}}

\begin{document}

\begin{abstract}
We construct a klt del Pezzo surface in characteristic three violating the Kawamata-Viehweg vanishing theorem. As a consequence we show that there exists a Kawamata log terminal threefold singularity which is not Cohen-Macaulay in characteristic three.
\end{abstract}

\maketitle

\tableofcontents

\section{Introduction}
In characteristic zero one of the main technical tool used to establish the Minimal Model Program (MMP for short) is the Kawamata-Viehweg vanishing theorem.
Unfortunately, vanishing theorems are known to fail in general for varieties defined over fields of positive characteristic and a great amount of work has been done to construct examples of pathological varieties violating them  and to study their geography (see for example \cite{Ray78}, \cite{Eke88}, \cite{Muk13} and \cite{dCF15}). 

In this context, varieties of Fano type over perfect fields of characteristic $p>0$ violating Kodaira-type vanishing theorems seem rather rare and in fact are conjectured to exist only for small primes in each dimension. 
As far as the author knows, only two classes of varieties of Fano type violating Kawamata-Viehweg vanishing have been constructed and both are in characteristic two (see \cite{LR97} and \cite{Kov18} for a six-dimensional smooth example and \cite{CT} for a two-dimensional one with klt singularities).
Moreover, the only vanishing statement supporting the above conjecture in every dimension is \cite[Theorem 1.9]{PW}, where it is proven that the classical Kodaira vanishing theorem holds for the first cohomology group on a klt Fano variety $X$ if $p>2\dim X$.

In the case of surfaces, however, the situation is much better understood: it is known that the Kawamata-Viehweg vanishing theorem holds for smooth del Pezzo surfaces over an algebraically closed field (see \cite[Appendix A]{CT18}), for regular del Pezzo surfaces over a (possibly imperfect) field of characteristic $p>3$ (see \cite[Theorem 1.1]{Das}) and for surfaces of del Pezzo type over an algebraically closed field of sufficiently large characteristic (see \cite[Theorem 1.2]{CTW17}), although we do not know an explicit lower bound on the characteristic.

In this article, inspired by an example of Keel and M$^\text{c}$Kernan (see \cite[Section 9]{KM99} and \cite[Section 4]{CT}), we construct a klt del Pezzo surface violating the Kawamata-Viehweg vanishing in characteristic three:
\begin{thm}[See Theorem \ref{failurekod3}]\label{thm1} 
Let $k$ be a field of characteristic three.
Then there exists a projective normal $k$-surface $T$ such that
\begin{enumerate}
\item $T$ has klt singularities and $-K_T$ is ample;
\item $\rho(T)=1$;
\item there exists an ample $\Q$-Cartier Weil divisor $A$ on $T$ such that $H^1(T,\mathcal{O}_T(-A)) \neq 0.$
\end{enumerate}
\end{thm}
\indent One important application of the Kawamata-Viehweg vanishing theorem for the MMP in characteristic zero is the proof due to Elkik that klt singularities are Cohen-Macaulay and rational (see \cite{Elk81}). 
In characteristic $p>0$, due to the failure of vanishing theorems, general cohomological properties of klt singularities are still largely unknown but, according to a local-global principle, they are expected to be strictly related to vanishing theorems for varieties of Fano type.

In dimension three, the main result of \cite{HW} shows that klt threefold singularities are Cohen-Macaulay and rational in large characteristic and the main ingredient of their proof is the Kawamata-Viehweg vanishing for surfaces of del Pezzo type. 
As for low characteristic, in \cite[Theorem 1.3]{CT} the authors give an example of a klt not Cohen-Macaulay threefold in characteristic two.  
Using a generalized cone construction and Theorem \ref{thm1} we give the first example of a klt singularity which is not Cohen-Macaulay in characteristic three, answering a question of Hacon-Witaszek and Kov\'acs (\cite[Question 5.4]{HW} and \cite{Kov18}):
\begin{thm}[See Section \ref{s-kltnotCM}] \label{klt not CM}
Let $k$ be a field of characteristic three.
Then there exists a normal $k$-variety $X$ such that
\begin{enumerate}
\item $X$ is a $\Q$-factorial variety of dimension three;
\item $X$ is klt;
\item $X$ is not Cohen-Macaulay (and in particular the singularities of $X$ are not rational).
\end{enumerate}
\end{thm}
An open problem is whether the Kawamata-Viehweg vanishing holds for surfaces of del Pezzo type over perfect fields of characteristic $p \geq 5 $. Despite not being able to solve this problem, in the last section we present a Kodaira-type vanishing theorem for big and nef Cartier divisors on klt del Pezzo surfaces of characteristic $p \geq 5$, answering a question of Cascini and Tanaka (see \cite[Remark 3.2]{CT18}).  
\begin{thm}[See Theorem \ref{vanishing}]\label{vanishingCartier}
Let  $k$ be an algebraically closed field of characteristic $p \geq 5$.
Let $X$ be a klt del Pezzo surface over $k$ and let $A$ be a big and nef Cartier divisor on $X$.  Then,
\[H^1(X,\mathcal{O}_X(A))=0. \]
\end{thm}
\begin{remark}
After a preliminary draft of this paper was announced, B. Totaro in \cite{Tot17} constructed for every $p>0$ a smooth Fano variety of dimension $2p+1$ violating the Kodaira vanishing theorem. As a consequence he constructs a terminal not Cohen-Macaulay singularity of dimension $2p+2$. 
\end{remark}
\indent \textbf{Acknowledgements:} 
I would like to express my gratitude to my advisor Paolo Cascini for his constant support and his useful suggestions and to Hiromu Tanaka, since this project started from a conversation I had with him on Keel-M$^{\text{c}}$Kernan surfaces and for answering my questions. I would also like to thank Ivan Cheltsov, Mirko Mauri, Omprokash Das, Davide Cesare Veniani and Jakub Witaszek for reading an earlier draft and for fruitful conversations. This work was supported by the Engineering and Physical Sciences Research Council [EP/L015234/1].

\vspace{2mm}
\section{Preliminaries}
\subsection{Notation} 
\begin{enumerate}
\item Throughout this article, $k$ denotes a field of characteristic $p \geq 0$. If $k$ is a perfect field of characteristic $p>0$, we denote by $W_m(k)$ the ring of Witt vectors of length $m$. 
\item By \emph{variety} we mean an integral scheme which is separated and of finite type over $k$.  If $X$ is a normal variety, we denote by $K_X$ its canonical divisor class (for more details on how to define the canonical class on singular varieties, we refer to \cite[Definition 1.6]{Kol13}).
\item Given a scheme $X$ defined over a field $k$ of characteristic $p>0$, we denote by $F \colon X \rightarrow X$ the absolute Frobenius morphism and for any natural number $e>0$ we denote the $e$-th iterate of Frobenius by $F^e$.
\item We say that $(X,\Delta)$ is a \emph{log pair} if $X$ is a normal variety, $\Delta$ is an effective $\Q$-divisor and $K_X+\Delta$ is a $\Q$-Cartier $\Q$-divisor. 
If $k$ is perfect, we say the pair $(X, \Delta)$ is \emph{log smooth} if $X$ is smooth and $\Supp(\Delta)$ is a snc divisor. We refer to \cite{KM98} and \cite{Kol13} for the definition of the singularities appearing in the MMP (e.g. \emph{klt, plt, lc}). We say $f \colon Y \rightarrow X$ is a \emph{log resolution} if $f$ is a proper birational morphism, $Y$ is smooth, $\Exc(f)$ has pure codimension one and the pair $(Y, \Supp(f^{-1}\Delta + \Exc(f)))$ is log smooth.
\item We say that a normal projective surface $X$ is  a \emph{surface of del Pezzo type} if there exists an effective $\Q$-divisor $\Delta$ such that $(X, \Delta)$ is klt and $-(K_X+\Delta)$ is ample. The pair $(X, \Delta)$ is said to be a log del Pezzo pair. We say $X$ is a klt del Pezzo surface if the pair $(X,0)$ is log del Pezzo.
\item  Let $X$ be a normal variety. Given $D$ a Weil divisor, we define the reflexive sheaf $\mathcal{O}_X(D)$ which to an open subset $U \subset X$ associates $H^0(U,\mathcal{O}_X(D))= \left\{ \varphi \in K(X) \mid \text{div}(\varphi)|_{U} +D|_{U} \geq 0 \right\}.$
If $L$ is a reflexive sheaf of rank one on $X$, there exists a Weil divisor $D$ such that $L \simeq \mathcal{O}_X(D)$. We denote by $L^{[m]}$ the double dual of $L^{\otimes m}$, which is isomorphic to $\mathcal{O}_X(mD)$. 
By $\Cl(X)$ we denote the abelian group of $\Z$-Weil divisors modulo linear equivalence.
If $X$ is a proper variety over $k$,  we denote by $\rho(X):= \rho(X/\Spec \,k)$ the \emph{Picard number} of $X$.
\item For us, an $(A_n)$-$klt$ (resp. $(A_n)$-$canonical$) singularity is a klt (resp. canonical) surface singularity such that the exceptional locus of its minimal resolution is a string of $n$ smooth rational curves. 
\end{enumerate}

\subsection{Frobenius splitting}
We fix a perfect field $k$ of characteristic $p>0$. 
For the convenience of the reader we recall the definition of $F$-splitting.
\begin{definition}
Let $X$ be a variety over $k$. We say that $X$ is \emph{globally F-split} if for some $e >0$ the natural map
\begin{equation*}
\xymatrix{
\mathcal{O}_X \ar[r]^{F^e} & F^e_* \mathcal{O}_X 
}
\end{equation*}
splits as a homomorphism of $\mathcal{O}_X$-modules.
\end{definition}
\begin{remark}\label{splitt}
In the definition of $F$-splitting for a variety $X$ one can equivalently ask that for \emph{all} $e>0$ the $\mathcal{O}_X$-module homomorphism $\mathcal{O}_X \rightarrow  F^e_* \mathcal{O}_X $ splits. 
\end{remark}
Being globally $F$-split implies strong vanishing results for the cohomology of ample divisors on $X$. We will need the following result, which is a mild generalization of \cite[Theorem 1.2.9]{BK05} to $\Q$-Cartier Weil divisors on surfaces:

\begin{proposition}\label{KodairaforFsplit}
Let $X$ be a normal projective surface over ~$k$. If $X$ is globally $F$-split, then for any ample $\Q$-Cartier Weil divisor $A$
\[H^1(X, \mathcal{O}_X(-A))=0. \]
\end{proposition}

\begin{proof}
By Remark \ref{splitt} we know that for any large $g \gg 0$ there exists a splitting:
\[\mathcal{O}_X \rightarrow F^g_* \mathcal{O}_X \rightarrow \mathcal{O}_X. \]
Restricting to the regular locus $U$ and tensoring by $\mathcal{O}_U(-A)$ we have the following splitting:
\[\mathcal{O}_U(-A) \rightarrow F^g_* \mathcal{O}_U(-p^g A) \rightarrow \mathcal{O}_U(-A). \]
Since $X$ is a normal variety and each sheaf in the sequence is reflexive we deduce that the splitting holds on the whole $X$:
\[\mathcal{O}_X(-A) \rightarrow F^g_* \mathcal{O}_X(-p^g A) \rightarrow \mathcal{O}_X(-A). \]
Passing to cohomology, we have an injection:
\[H^1(X,\mathcal{O}_X(-A) ) \hookrightarrow H^1(X, \mathcal{O}_X(-p^gA)).\]
Let $m$ be the Cartier index of $A$ and let us write $m=p^f h$ where $\gcd (p,h)=1$. Then for large enough and sufficiently divisible $e$ we have that $m$ divides $p^f(p^e-1)$ and thus
\[p^f(p^e-1)A \text{ is Cartier}. \]
Now consider $g=f+e$. We have
\[-p^g A= -(p^{e+f}-p^f)A-p^f A.\]
Since $X$ is a normal surface, any divisorial sheaf is Cohen-Macaulay and thus we can apply Serre duality (see \cite[Theorem 5.71]{KM98}) to deduce
\[H^1(X, \mathcal{O}_X(-p^gA))\simeq H^{1}(X,\mathcal{H}om_{\mathcal{O}_X}(\mathcal{O}_X(-p^gA), \mathcal{O}_X(K_X)))^{*}. \]
Being $X$ normal, we have $\mathcal{H}om_{\mathcal{O}_X}(\mathcal{O}_X(-p^gA), \mathcal{O}_X(K_X))\simeq \mathcal{O}_X(K_X+p^gA)$   since both sheaves are reflexive and they are isomorphic on the regular locus $U$. Therefore
\begin{align*}
H^1(X, \mathcal{O}_X(-p^gA)) &\simeq H^{1}(X, \mathcal{O}_X(K_X+p^gA))^*   \\
&=H^{1}(X, \mathcal{O}_X(K_X+p^fA) \otimes \mathcal{O}_X((p^{e+f}-p^f)A))^*.
\end{align*}
By choosing $e$ sufficiently large and divisible we conclude that the last cohomology group vanishes by the Serre vanishing criterion for ample line bundles.
\end{proof}

\subsection{A cone construction for Weil $\Q$-Cartier divisors}
\label{ss-cone-const}

For the theory of cones on polarised algebraic varieties $(X, L)$ where $X$ is a projective normal variety and $L$ is an ample Cartier divisor, we refer to \cite[Chapter 3]{Kol13}.
In Section \ref{s-kltnotCM} we need to deal with a generalisation, originally due to Demazure \cite{Dem88}, of the cone construction to the case of ample $\Q$-Cartier Weil divisors.
We thus extend some of the results explained by Koll\'ar to this setting.  

Let $(X,\Delta)$ be a log pair of dimension $n$ where $X$ is a $k$-projective variety and let $L$ be an ample $\Q$-Cartier Weil divisorial sheaf on $X$. 
The variety
\[C_a(X,L):=\Spec_k \sum_{m \geq 0} H^0 (X,L^{[m]})\] 
is the cone over $X$ induced by $L$.
The closed point defined by the ideal $\sum_{m \geq 1} H^0(X,L^{[m]})$ is called the \emph{vertex} of the cone and we denote it by $v$.
Over $X$ we consider the affine morphism:
\[\pi \colon BC_a(X,L) := \Spec_X \sum_{m \geq 0} L^{[m]} \rightarrow X.\]
The morphism $\pi$ comes with a natural section $X^-$ defined by the vanishing of the ideal sheaf $\sum_{m \geq 1} L^{[m]}$. 
The open subset of $BC_a(X,L)$
\[\Spec_X \sum_{m \in \Z} L^{[m]}= BC_a(X,L) \setminus X^-, \]
is isomorphic to $C_a^*(X,L):=C_a(X,L) \setminus v$. 
We have the following diagram:
\begin{equation} \label{diagramcone}
\xymatrix{
	 BC_a(X,L) \ar[r]^{\qquad \pi} \ar[d]_{f} & X  \\
	 C_a(X,L) &.
}
\end{equation}

The birational morphism $f$ contracts exactly the section $X^{-}$ of $\pi$ with anti-ample $\Q$-Cartier divisorial sheaf $\MO_{X^-}(X^-) \simeq L^{\vee}$.
Given a $\Q$-Cartier $\Q$-divisor $D$ on $X$, we can construct a $\Q$-divisor $D_{C_a(X,L)}=f_*\pi^* D$ on $C_a(X,L)$.

The following result (originally due to \cite{Wat81}) describes the divisor class group of the cone and the condition under which the log canonical class is $\Q$-Cartier.
For sake of completess, we present a proof.

\begin{proposition}[cf. {\cite[Proposition 3.14]{Kol13}}]\label{p-Qcart-div}
With the same notation as in diagram (\ref{diagramcone}), we have
\begin{enumerate}
\item $\Pic(C_a(X,L))=0 $,
\item $\Cl(C_a(X,L))=\Cl(X)/ \langle L \rangle$,
\item $m(K_{C_a(X,L)}+\Delta_{C_a(X,L)})$ is Cartier if and only if there exists $r \in \Q$ such that $\mathcal{O}_X(m(K_X+\Delta)) \simeq L^{[rm]}$. \label{canonicalclass}
Moreover,
\begin{equation} \label{discrep}
K_{BC_a(X,L)}+\pi^*\Delta+(1+r)X^-=f^*(K_{C_a(X,L)}+\Delta_{C_a(X,L)}).
\end{equation} 
\end{enumerate}
\end{proposition}
\begin{proof}
Let $Z$ be the locus where $L$ is not a Cartier divisor and denote the open set $V:= X \setminus Z$. 
Since $X$ is normal, $Z$ has codimension at least two and thus 
$\Cl(V) \simeq \Cl(X)$ and $\Pic(V) \simeq \Pic(X)$.
By construction, $\pi \colon \pi^{-1}(V) \to V$ is an $\mathbb{A}^1$-bundle and thus $\Cl(\pi^{-1}(V)) \simeq \Cl(V)$ and  $\Pic (\pi^{-1}(V)) \simeq \Pic(V)$.
Since $\pi$ is equi-dimensional, we conclude that the codimension of $\pi^{-1}(Z)$ is at least two and thus 
\[ \Cl(BC_a(X,L)) \simeq \Cl(\pi^{-1}(V)) \simeq \Cl(V) \simeq \Cl(X), \]
and analogously we have $\Pic(BC_a(X,L)) \simeq \Pic(X)$. 

We prove (1). Let $M$ be a line bundle on $C_a(X,L)$. 
Then $f^*M$ is trivial on $X^{-}$, thus concluding that $M$ is trivial on $C_a(X,L)$.

We prove (2). 
Since $v$ has codimension at least two in $C_a(X, L)$, we have $\Cl(C_a(X,L)) \simeq \Cl(C_a^*(X,L))$.  
Thus we have
\[ \mathbb{Z}[X^{-}] \to \Cl(BC_a(X,L)) \to \Cl(C_a^*(X,L)) \to 0. \]
Since $\MO_{X^{-}}(X^-) \simeq L^{\vee}$, we conclude (2).

We prove (3). 
We have $\pi^{-1}(V) \setminus X^{-} \to V$ is a $\mathbb{G}_m$-bundle.
Thus there is a natural linear equivalence $K_{\pi^{-1}(V)} +X^{-}|_{\pi^{-1}(V)} \sim \pi^*K_V$ which extends to a natural linear equivalence
$ K_{BC_a(X,L)} +X^{-} \sim \pi^*K_X $.
Thus $K_{BC_a(X,L)} +\pi^*\Delta +X^{-} \sim \pi^*(K_X+ \Delta).$
By assertion (2), the divisor $mK_{C_a(X,L)} + m\Delta_{C_a(X,L)}$ is Cartier if and only if $mK_{C_a^*(X,L)} + m\Delta_{C_a^*(X,L)}$ is linearly equivalent to zero.
This is equivalent to $m(K_X+ \Delta) \sim rm L$ for some $r \in \mathbb{Q}$.
As for the last equality, we have
\[ K_{BC_a(X,L)}+\pi^*\Delta+(1+a)X^-=f^*(K_{C_a(X,L)}+\Delta_{C_a(X,L)}), \]
for some $a$.
By restricting to $X^-$ we have
\begin{small}
\[ 0 \sim_{\Q} (K_{BC_a(X,L)}+\pi^*\Delta+(1+a)X^{-})|_{X^{-}} \sim_{\Q} (K_X + \Delta + a X^{-})|_{X^{-}} \sim_{\Q}  rL-aL, \]
\end{small}
thus concluding.
\end{proof}

From the point of view of the singularities of the MMP we have the following
\begin{proposition}[cf. {\cite[Lemma 3.1]{Kol13}}]\label{plt}
With the same notation as in diagram (\ref{diagramcone}), let us assume that $K_X + \Delta \sim_{\Q} rL$ for some $r \in \mathbb{Q}$. 
\begin{enumerate}
\item If $L$ is Cartier, then the pair $(C_a(X,L), \Delta_{C_a(X,L)})$ is terminal (resp. klt) if and only if the pair $(X,\Delta)$ is terminal (resp. klt) and $r < -1$ (resp. $r <0$).
\item If $L$ is Cartier, then the pair $(C_a(X,L), \Delta_{C_a(X,L)})$ is dlt if the pair $(X,\Delta)$ is dlt and $r <0$.
\item If $X$ is $\Q$-factorial and $\car(k)=0$, then $(C_a(X,L), \Delta_{C_a(X,L)})$ is klt if and only if $(X, \Delta)$ is klt and $r <0$.
\end{enumerate}
\end{proposition}
\begin{proof}
Assertions (1) and (2) are proven in \cite[Lemma 3.1]{Kol13}.
We prove assertion (3).
Since $X$ is $\Q$-factorial we have that $BC_a(X,L)$ is $\Q$-factorial.
Since $(X, \Delta)$ is klt, we conclude by inversion of adjunction (see \cite[Theorem 5.50]{KM98}) that the pair $(BC_a(X,L), \pi^*\Delta + X^{-})$ is plt and the unique plt centre is $X^{-}$.
Since $r<0$, we have $(BC_a(X,L), \pi^*\Delta + (1+r )X^{-})$ is klt and thus by Proposition \ref{p-Qcart-div} we conclude the proof.
\end{proof}

We will be interested in understanding whether the singularity at the vertex of the cone is Cohen-Macaulay or not.
For this reason, we show that the local cohomology at the vertex of the cone is controlled by the cohomology groups of $L$ and its multiples:

\begin{proposition} \label{local cohomology} For $i \geq 2$,
\[H^i_v(C_a(X,L),\mathcal{O}_{C_a(X,L)}) \simeq  \bigoplus_{m \in \Z} H^{i-1}(X,L^{[m]}).\]
\end{proposition}
\begin{proof}
Since $C_a(X,L)$ is affine, we have that the cohomology groups $H^i(C_a(X,L), \mathcal{O}_{C_a(X,L)})$ vanish for $i \geq 1$. Thus, by the long exact sequence in local cohomology (see \cite[Chapter III, ex. 2.3]{Har77}) we deduce
\[H^i_v(C_a(X,L),\mathcal{O}_{C_a(X,L)}) \simeq H^{i-1}(U, \mathcal{O}_U) \text{ for } i \geq 2,\]
where $U:= C_{a}^*(X,L) \simeq \sum_{m \in \mathbb{Z} }\Spec_X L^{[m]}.$
Since $\pi$ is an affine morphism, we have
\begin{small}
\begin{equation*}
H^{i-1}(U, \mathcal{O}_U)=H^{i-1}(X,\pi_* \mathcal{O}_U)=
H^{i-1}(X,\bigoplus_{m \in \Z} L^{[m]})= \bigoplus_{m \in \Z} H^{i-1}(X,L^{[m]}),
\end{equation*}
\end{small}
thus concluding.
\end{proof}
\section{A Keel-M$^{\text{c}}$Kernan surface in characteristic three} \label{s-logdPKVV}
In this section and in the following we fix $k$ to be a field of characteristic three. We prove Theorem \ref{thm1} by constructing a klt del Pezzo surface of Picard rank one not satisfying the Kawamata-Viehweg vanishing theorem.

\subsection{Construction}
In \cite[Section 9]{KM99}, the authors construct a family of klt del Pezzo surfaces in characteristic two violating the Bogomolov bound on the number of singular points. 
In \cite{CT} it was noted that their example gives various counterexamples to the Kawamata-Viehweg vanishing theorem. We adapt their construction to the case of characteristic three. 

Let us consider the smooth rational curve $C$ inside $\mathbb{P}^1_x \times \mathbb{P}^1_y$ defined by the equation:
\[C:= \left\{([x_1:x_2],[y_1:y_2]) \mid x_2y_1^3-x_1y_2^3 =0 \right\}. \]
We denote by $\pi_x \colon \mathbb{P}^1_x \times \mathbb{P}^1_y \rightarrow \mathbb{P}^1_x$ the natural projection onto the first coordinate and we say $F_p=\pi_x^{-1}(p)$ for $p \in \mathbb{P}^1_x$ is the vertical fibre over $p$. \\
\indent The main property of $C$ we are interested in is that, since the characteristic of the field $k$ is three, the morphism $\pi_x|_C \colon C \rightarrow \mathbb{P}^1_x$ is the geometric Frobenius morphism. Geometrically, the curve $C$ has the following ``funny'' property: every vertical fibre $F_p$ is a triple tangent to $C$. \\
\indent Fix a closed $k$-point $p_1$ on $C$ and consider the vertical fiber $F_1$ passing through this point. Since such a fiber is a triple tangent to $C$ at the point $p_1$ we perform three successive blow-ups to separate $C$ from $F_1$. The order of the blow-ups is as follows: at each step we blow up the intersection point of the strict transform of $F_1$ and the strict transform of $C$. After these birational modifications the strict transforms of $C$ and $F_1$ (which, by abuse of notation, are denoted by the same letter) and the exceptional divisors $E_1,G_1,H_1$ are in the following configuration: 
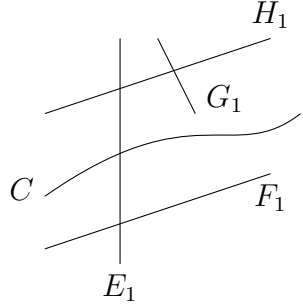
\begin{center}
\begin{tikzpicture}
\draw (0,1) -- (0,4);
\draw (-1,1.2) -- (2,2.2);
\draw (-1,1.9) .. controls (1,3.3) and (1.5,2.3) .. (2.4,3);
\draw (-1,3) -- (2, 4);
\draw (1,3) -- (0.5,4);
\node [below] at (0,1) {$E_1$};
\node [below] at (2,2.2) {$F_1$};
\node [above] at (2,4) {$H_1$};
\node [right] at (1,3.2) {$G_1$};
\node [left] at (-1,2) {$C$};
\end{tikzpicture}
\end{center}
where all the curves are smooth and rational with the following intersection numbers:
\begin{align*}
H_1^2=-2, G_1^2=-2,  F_1^2=-3, E_1^2=-1, \\
C \cdot E_1=1, E_1 \cdot F_1=1, E_1 \cdot H_1=1, H_1 \cdot G_1=1.
\end{align*}
Note that the self-intersection of $C$ has dropped by three. Performing the same operation with other two $k$-points $p_2, p_3$ on the curve $C$ we construct a birational morphism $f \colon S \rightarrow \mathbb{P}^1_x \times \mathbb{P}^1_y$ where the strict transform of $C$ has become a $(-3)$-curve. Over each point $p_i$ we have the exceptional curves $H_i, G_i, E_i$ and the strict transform of the fibre $F_i$ in the same configuration as the one described above for $p_1$. \\
\indent On $S$ there are the $(-3)$-curves $F_1, F_2, F_3$ and three disjoint strings of two $(-2)$-curves formed by $H_i$ and $G_i$ for $i=1, \dots, 3$. Let $\psi \colon S \rightarrow T$ be the birational contraction of the curves $F_i, H_i, G_i$ for $i \in \left\{1,2,3 \right\}$ and $C$. We can construct $\psi$ by running a suitable log MMP (see \cite{Tan14}) for the log pair \[(S,\Delta:=\sum_{i=1}^3 \frac{2}{3}F_i + \sum_{i=1}^3 \frac{1}{2}(H_i+G_i)+\frac{2}{3}C),\]
which at each step of the MMP we contract exactly one of the curves appearing in $\Supp(\Delta)$. We denote, with a slight abuse of notation, the pushforward of a divisor $D$ via $\psi$ with the same letter $D$.\\
\indent On $T$ we have the following configuration of curves and singular points:\\
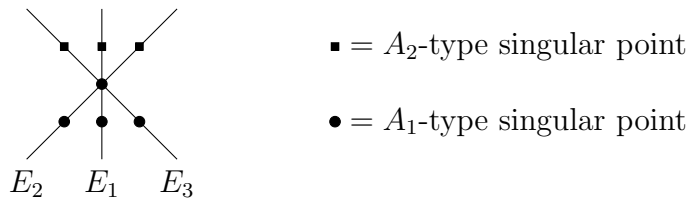
\begin{center}
\begin{tikzpicture}
\draw (0,1) -- (0,3);
\draw (-1,1) -- (1, 3);
\draw (1,1) -- (-1, 3);
\node [below] at (0,1) {$E_1$};
\node [below] at (-1,1) {$E_2$};
\node [below] at (1,1) {$E_3$};
\filldraw[fill=black, draw=black] (-0.05,2.45) rectangle (0.05,2.55);
\filldraw[fill=black, draw=black] (-0.55,2.45) rectangle (-0.45,2.55);
\filldraw[fill=black, draw=black] (+0.55,2.45) rectangle (0.45,2.55);
\draw [fill=black] (0,2) circle [radius=0.07];
\draw [fill=black] (0.5,1.5) circle [radius=0.07];
\draw [fill=black] (0,1.5) circle [radius=0.07];
\draw [fill=black] (-0.5,1.5) circle [radius=0.07];
\draw [fill=black] (-0.5,2.5);
\filldraw[fill=black, draw=black] (3.05,2.45) rectangle (3.15,2.55);
\node [right] at (3.15,2.50) {$=A_2$-canonical singularity};
\draw [fill=black] (3.10,1.5) circle [radius=0.07];
\node [right] at (3.15,1.5) {$=A_1$-klt singularity};
\end{tikzpicture}
\end{center}

\begin{remark}\label{singularpoints}
The singularity at the points of type $A_2$-canonical (resp. $A_1$-klt) is formally isomorphic to the quotient of $\mathbb{A}^2_k$ by the action of the group scheme $\mu_3$ with weights $(1,2)$ (resp. $(1,1)$).
\end{remark}

The following proposition justifies why this surface is a generalization of Keel-M$^{\text{c}}$Kernan's example in characteristic three:

\begin{proposition} \label{p-T-is-dP}
The surface $T$ is a klt del Pezzo surface of Picard rank one.
Moreover, $-K_T \equiv E_1$. 
\end{proposition}
\begin{proof}
It is straightforward to see that $\rho(T)=1$. Since we contract only cycle of $(-2)$-curves and $(-3)$-curves, $T$ has klt singularities. We are only left to show that $-K_T$ is an ample divisor. By an explicit computation we have
\begin{equation} \label{canonical}
\psi^* K_T=K_S+\sum_{i=1}^3 \frac{1}{3} F_i +\frac{1}{3}C.
\end{equation}
Since $\rho(T)=1$ and $h^1(T, \MO_T)=0$ it is enough to prove that the anticanonical divisor has the same intersection as $E_1$ when intersected with any effective curve. \\
\indent Let $F_p$ be the fibre of the the map $\pi_x \circ f \colon S \rightarrow \mathbb{P}^1_x$ over a general point $p \in \mathbb{P}^1_x$. By the projection formula we have:
\[ -K_T \cdot \psi_* F_p =-\psi^*K_T \cdot F_p=-K_S \cdot F_p-\frac{1}{3} C \cdot F_p=1. \qedhere\]
Again by the projection formula, we also have
\[E_1 \cdot \psi_* F_p = \psi^* E_1 \cdot F_p = \frac{1}{3} C \cdot F_p=1.  \]
\end{proof}
\begin{remark}
It is possible to perform a similar construction for higher characteristic, but the resulting surface will have ample canonical divisor class. 
\end{remark}
\begin{remark}\label{notlift}
In \cite[Section 9]{KM99} the authors prove the Bogomolov bound: a klt del Pezzo surface of Picard rank one over an algebraically closed field of characteristic zero has at most six singular points. The bound was later improved to four singular points in characteristic zero in \cite{Bel09}. The surface $T$ has seven singular points and thus shows that the Bogomolov bound does not hold in characteristic three. 
It is an open question whether the Bogomolov bound holds for large characteristic.
\end{remark}
We show that there are no anticanonical sections on $T$:
\begin{proposition} \label{-Kiseffective}
$H^0(T, \mathcal{O}_T(-K_T)) = 0$.
\end{proposition}
\begin{proof}
By formula (\ref{canonical})
we have 
\[H^0(T, \mathcal{O}_T(-K_T))= H^0(S, \mathcal{O}_S(-K_S- \sum_{i=1}^3 F_i -C)).\]
A direct computation shows
\begin{small}
\[-K_S-\sum_{i=1}^3 F_i-C \sim f^*({-K_{\mathbb{P}^1_x \times \mathbb{P}^1_y}}-\sum_{i=1}^3 F_i-C)+ \sum_{i=1}^3 (G_i+2H_i+3E_i). \]
\end{small}
Therefore
\begin{align*}
H^0(T, \mathcal{O}_T(-K_T))&=H^0(\mathbb{P}^1_x \times \mathbb{P}^1_y, \mathcal{O}(-K_{\mathbb{P}^1_x \times \mathbb{P}^1_y}-\sum_{i=1}^3 F_i-C)) \\
&=H^0(\mathbb{P}^1_x \times \mathbb{P}^1_y, \mathcal{O}(-2,-1))=0.
\end{align*}
\end{proof}
\subsection{Failure of the Kawamata-Viehweg vanishing theorem}
We show that the Kawamata-Viehweg vanishing theorem fails on the surface $T$.

We consider the following ample $\Q$-Cartier Weil divisor 
\[ A:=E_2+E_3-E_1. \]
\begin{thm}\label{failurekod3}
The Kawamata-Viehweg vanishing theorem fails for the Weil divisor $A$; \emph{i.e.}
\[H^1(T,\mathcal{O}_T(-A)) \neq 0.\]
\end{thm}
\begin{proof}
The strategy is to pull-back the divisor to the minimal resolution $S$ and compute there the cohomology groups. Let us consider the pull-back of $A$ to $S$ as a $\Q$-divisor:
\begin{small}
\[-\psi^*A=E_1+\frac{1}{3}F_1+\frac{2}{3} H_1+\frac{1}{3}G_1- E_2-\frac{1}{3}F_2-\frac{2}{3} H_2-\frac{1}{3}G_2-E_3-\frac{1}{3}F_3-\frac{2}{3}H_3-\frac{1}{3}G_3-\frac{1}{3}C;\]
\end{small}
%\begin{small}
%\[-\psi^*A=E_1+\frac{1}{3}F_1+\frac{2}{3} H_1+\frac{1}{3}G_1+\frac{1}{3}C- E_2-\frac{1}{3}F_2-\frac{2}{3} H_2-\frac{1}{3}G_2-\frac{1}{3}C-E_3-\frac{1}{3}F_3-\frac{2}{3}H_3-\frac{1}{3}G_3-\frac{1}{3}C;\]
%\end{small}
thus
\[ \lfloor{-\psi^*A}\rfloor=E_1-E_2-F_2-H_2-G_2-E_3-F_3-H_3-G_3-C. \]
We have
\[ \psi_* \mathcal{O}_S(\lfloor{-\psi^*A}\rfloor)=\mathcal{O}_T(-A),\]
and we compute the cohomology group using the Leray spectral sequence
\begin{equation} \label{leray}
E^{i,j}_2= H^j(T, R^i \psi_*\mathcal{O}_S(\lfloor{-\psi^*A}\rfloor)) \Rightarrow H^{i+j}(S,\mathcal{O}_S(\lfloor{-\psi^*A}\rfloor)). 
\end{equation}

We show that $R^i\psi_* \mathcal{O}_S(\lfloor{-\psi^*A}\rfloor)=0$ for $i>0$.
By the Kawamata-Viehweg vanishing theorem for birational morphism between surfaces (see \cite[Theorem 10.4]{Kol13}) we just need to check that $\lfloor{-\psi^*A}\rfloor$ is $\psi$-nef: 
\begin{align*}
&\lfloor{-\psi^*A}\rfloor \cdot C=2, \\
& \lfloor{-\psi^*A}\rfloor \cdot F_1 =1,  \lfloor{-\psi^*A}\rfloor \cdot H_1=1, \lfloor{-\psi^*A}\rfloor \cdot G_1=0, \\
& \lfloor{-\psi^*A}\rfloor \cdot F_i=2, \lfloor{-\psi^*A}\rfloor \cdot H_i= 0, \lfloor{-\psi^*A}\rfloor \cdot G_i=1 \text{ for } i=2,3.
\end{align*}
Therefore the Leray spectral sequence (\ref{leray}) degenerates at the $E_2$-page and we have for all $i \geq 0$:
\[H^i(T, \mathcal{O}_T(-A)) \simeq H^i(S, \mathcal{O}_S(\lfloor{-\psi^*A}\rfloor)). \]
By a direct computation we have
\[ K_S \cdot \lfloor{-\psi^*A}\rfloor=-2 \text{ and } \lfloor{-\psi^*A}\rfloor^2=-6. \]
Therefore, by the Riemann-Roch theorem on $S$, we deduce 
\[ \chi(T, \mathcal{O}_T(-A))=\chi(S, \mathcal{O}_S(\lfloor{-\psi^*A}\rfloor))=-1, \]
which implies $h^1(T,\mathcal{O}_T(-A)) \neq 0$.
\end{proof}

We now conclude that the surface $T$ is neither $F$-split or admits a log resolution lifting to the second Witt vectors, giving thus a generalization of \cite[Theorem 1.3]{CTW17} to characteristic three. 
\begin{cor}
Over any perfect field $k$ of characteristic $p=3$ there exists a log del Pezzo surface $T$ which is not globally $F$-split and such that for any log resolution $\mu \colon S \rightarrow T$ the log smooth pair $(S,\Exc(\mu))$ does not lift to $W_2(k)$. 
\end{cor} 

\begin{proof}
The surface $T$ constructed above is not globally $F$-split by Proposition \ref{KodairaforFsplit} and Theorem \ref{failurekod3}. 
By Proposition \ref{failurekod3} and Serre duality we have
\[H^1(T, \mathcal{O}_T(K_T+A)) \neq 0. \]
If the pair $(S, \Exc(\mu))$ lifted to $W_2(k)$, we could apply \cite[Lemma 6.1]{CTW17} to the $\Z$-divisor $D:=K_T+A$, thus getting a contradiction with the non-vanishing above.
\end{proof}
\section{A klt threefold singularity not CM in characteristic three} \label{s-kltnotCM}

In this section we construct an example of a klt threefold singularities in characteristic three which is not Cohen-Macaulay.

With the same notation as in Subsection \ref{ss-cone-const}, let us consider the cone over the klt del Pezzo surface $T$ constructed in Section \ref{s-logdPKVV}:
\[X:=C_a(T,\mathcal{O}_T(A)),\]
where $A$ is the $\Q$-Cartier Weil divisor of Theorem \ref{failurekod3}.

We now prove that $X$ has klt singularities. 
If we were working over a field of characteristic zero we would conclude immediately by Proposition \ref{plt}.
However, since we are working in positive characteristic, we need to further study the singularities of $X$ to conclude it is klt.
We start by studying the singularities of its partial resolution:
\[Y:=  \Spec_T \sum_{m \geq 0}  \mathcal{O}_T(mA) \xrightarrow{f} X. \]
The exceptional locus of the birational morphism $f$ is the prime divisor $E$, which is isomorphic to $T$. We denote by $\pi$ the natural affine map $\pi \colon Y \rightarrow T$. We thus have the following diagram:
\begin{equation*}
\xymatrix{
	 Y \ar[r]^{\pi} \ar[d]_{f} & T  \\
	 X & .
}
\end{equation*}

\begin{proposition} \label{p-first-step}
The variety $Y$ is a $\Q$-factorial threefold with isolated singularities and the pair $(Y,E)$ is toroidal (hence  log canonical).
\end{proposition}

\begin{proof}
To check that $Y$ is $\Q$-factorial it is sufficient to work in an analytic neighbourhood of the singular locus by \cite[(24.E)]{Mat80}. The same is true to compute the discrepancies. Thus we can reduce to study the preimage $\pi^{-1}(U) \subset Y$ of an analytic neighbourhood $U$ of the singular points of $T$ because outside the preimage of those points the pair $(Y,E)$ is log smooth.\\
\indent As explained in Remark \ref{singularpoints}, there are two different types of singular points in $T$. We show the result is true for the $A_2$-type singular points; for the $A_1$-type singular points the computation is similar. \\ 
\indent Let us consider a singular point $p \in T$, which is formally isomorphic to the quotient of $\mathbb{A}^2_{u,v}$ by the group $\mu_3$ with weight $(1,2)$. In local coordinates,
\[\mathbb{A}^2_{u,v} \git \mu_3 =\Spec_k k[u^3, v^3, uv] \simeq \Spec_k k[x,y,z]/(z^3-xy), \]
and the Weil divisorial sheaf $\mathcal{O}_T(A)$ is isomorphic to the Weil divisorial ideal $D:=(x,z)$. In this case we have
\begin{align*}
Y &= \Spec_{\mathbb{A}^2 \git \mu_3} \sum_{m \geq 0} \mathcal{O}_T(mD) \\
& \simeq \Spec_k \frac{k[x,y,z,a,b,c,d]}{(z^3-xy, a^2-cx, ab-cz, a^3-dx^2, ac-dx, b^3-dy, bc-dz)}.
\end{align*}
The fibration $f$ is the natural morphism associated to the $k$-algebra homomorphism
\begin{small}
\[\frac{k[x,y,z]}{(z^3-xy)} \rightarrow  \frac{k[x,y,z,a,b,c,d]}{(z^3-xy, a^2-cx, ab-cz, a^3-dx^2, ac-dx, b^3-dy, bc-dz)},\]
\end{small}
and the section $E$ is the subvariety defined by the ideal $(a,b,c,d)$. \\
\indent A more conceptual way to understand the $\mathbb{A}^1$-fibration $Y \rightarrow T$ and its singularities is to see it locally as a quotient of the trivial $\mathbb{A}^1$-bundle over $\mathbb{A}^2$. Let us consider the line bundle 
\[\mathbb{L}:= \Spec_{\mathbb{A}^2} \sum_{k \geq 0} (u)^k \simeq \Spec_{k} k[u,v,s] \]
together with the section $S=(s=0)$. We have a natural action of $\mu_3$ on $\mathbb{L}$ of weight $(1,2,1)$ and we can construct the quotient
\[p \colon \mathbb{L} \rightarrow \mathbb{L} \git \mu_3.\]
A direct computation shows that the quotient pair $(\mathbb{L} \git \mu_3, p(S))$ is isomorphic to $(Y,E)$. With this description, we deduce that $Y$ is a $\Q$-factorial variety by \cite[Lemma 5.16]{KM98} and that the singularities of $Y$ are isolated. \\
\indent Moreover we have shown that, near the preimage via $f$ of the singular points of $T$, the pair $(Y,E)$ is toroidal and thus by \cite[Proposition 11.4.24]{CLS11} we conclude it has log canonical singularities.
\end{proof}

%To study discrepancies on $X$ we need the following:
%\begin{lem} \label{qgorenstein}
%The threefold $X$ is $\Q$-factorial.
%\end{lem}

%\begin{proof}
%Tensoring by $\Q$ the short exact sequence in Proposition \ref{divisor class group}, we have the short exact sequence of $\Q$-vector spaces
%\[0 \rightarrow \Q \rightarrow \Cl(T)_\Q \rightarrow \Cl(X)_\Q \rightarrow 0. \] 
%Since $\rho(T)=1$ and $H^1(T, \mathcal{O}_T)=0$, we have $\Cl(T)_{\Q} \simeq \Q$. Therefore $\Cl(X)_\Q =0$, concluding the proof. 
%\end{proof}

\begin{thm}\label{t-klt-notCM}
The variety $X$ has $\Q$-factorial klt singularities and it is not Cohen-Macaulay. 
\end{thm}
\begin{proof}
By Proposition \ref{local cohomology} and Theorem \ref{failurekod3} we deduce
\[H^2_v(X, \mathcal{O}_X) \simeq \sum_{m \in \Z} H^1(T,\mathcal{O}_T(mA)) \neq 0,\]
thus proving $X$ is not Cohen-Macaulay. We are left to check it is klt. 
We have $-K_T \sim_{\Q} A$ by Proposition \ref{p-T-is-dP}.
Thus by Proposition \ref{p-Qcart-div}, $K_X$ is $\Q$-Cartier and $K_Y \sim_{\Q} f^*K_X$.

By Proposition \ref{p-first-step}, $Y$ is a $\Q$-factorial variety and the pair $(Y,E)$ is toroidal.
Therefore by \cite[Proposition 11.4.24]{CLS11} $Y$ is klt and since $X$ is crepant to $Y$ we conclude that $X$ has klt singularities.
\end{proof}
\section{Kodaira-type vanishing for klt del Pezzo surfaces}
The aim of this section is to collect some Kodaira-type vanishing results for big and nef line bundles on klt del Pezzo surfaces for arbitrary $p>0$ and to prove Theorem \ref{vanishingCartier}. In particular, we answer a question of Cascini and Tanaka in the case of characteristic $p \geq 5$ (see \cite[Remark 3.2]{CT18})).

We start by discussing the case of klt del Pezzo surfaces with few non-canonical singular points. 
Let us recall the following result.
\begin{lem}\label{l-rational-dP}
Let $k$ be an algebraically closed field of characteristic $p>0$.
Let $X$ be a surface of del Pezzo type over $k$. 
Then $X$ is a rational surface and $H^i(X, \MO_X)=0$ for $i>0$.
In particular, $\chi(X, \MO_X)=1$.
\end{lem}
\begin{proof}
The surface $X$ is rational by \cite[Fact 3.4 and Theorem 3.5]{Tan15}.
Let $Y$ be the minimal resolution. By the Kawamata-Viehweg vanishing theorem for birational morphism between surfaces, we have
$H^i(X, \MO_X) \simeq H^i(Y, \MO_Y)$.
Since $Y$ is a smooth rational surface we conclude.
\end{proof}
\begin{proposition} \label{just4points}
Let $k$ be an algebraically closed field of characteristic $p>0$.
Let $X$ be a klt del Pezzo surface over $k$ with at most three  non-canonical singular points. 
Suppose that all of them are formally isomorphic to the quotient of $\mathbb{A}^2_k$ by the action of a group scheme $\mu_m$ for some $m>0$. 
Then for any big and nef Cartier divisor $A$ we have
\[H^1(X,\mathcal{O}_X(A))=0.\]
\end{proposition}
\begin{proof}
Let us recall that given a Weil divisor $D$ on a surface $X$ with only $\mu_m$-quotient singularities we have
\[ \chi(X,\mathcal{O}_X(D))= \chi(X, \mathcal{O}_X) + \frac{1}{2}D \cdot (D-K_X) + \sum_{P \in \text{NotCart(D)}} c_P (D), \]
where $c_P(D)$ is a rational number depending on the type of singularity of the pair $(X,D)$ near $P$ (for more details see \cite{Rei87}). 
In \cite[Corollary 4.1]{PV07}, the authors prove that $c_P (-K_X) > -1$ (let us note that the assumption on the characteristic of the base field is unnecessary). 
Applying the Riemann-Roch formula we have
\begin{align*}
h^0(\mathcal{O}_X(A-K_X)) & \geq 1+ \frac{1}{2}(A^2-3A \cdot K_X + 2 K_X^2) + \sum _{P \in \text{NotCart}(K_X)} c_P(-K_X) \\
& \geq 3 + K_X^2 -3 >0.
\end{align*}
Since $h^1(X, \MO_X)=0$, we can apply \cite[Proposition 3.3]{CT} to show that $h^1(X,\mathcal{O}_X(K_X-A))=0$. 
Thus, by Serre duality we conclude $h^1(X, \mathcal{O}_X(A))=0$.
\end{proof}
\begin{remark} \label{vanishingT}
By \cite[Corollary 4.1]{PV07}, we have that $c_P(-K_X) \geq -\frac{1}{3}$ if the singularity $P$ is formally a quotient of $\mathbb{A}^2_k$ by $\mu_3$. 
Then the same proof of the previous Proposition shows that vanishing theorems for big and nef Cartier divisors hold on the surface $T$ we constructed in Section \ref{s-logdPKVV}.
This explains why we had to look for a $\Q$-Cartier Weil divisor violating the vanishing theorem.
\end{remark}
\begin{remark}
Let us note that the Kawamata-Viehweg vanishing theorem is not valid for del Pezzo surfaces with canonical singularities in characteristic two by \cite[Theorem 3.1]{CT}.
\end{remark}
\indent We recall a result of Koll\'{a}r (see \cite[Theorem II.6.2, Remark II.6.2.4 and Remark II.6.7.2]{Kol96}), which is a combination of Ekedahl's purely inseparable trick and Bend and Break techniques.
\begin{thm}\label{kollarbend}
Let $k$ be an algebraically closed field of characteristic $p>0$.
Let $X$ be a projective normal surface over $k$. Let $L$ be a big and nef Weil $\Q$-Cartier divisor on $X$ such that $H^1(X,L^{\vee})\neq 0$. Assume that $X$ is covered by a family of curves $\left\{D_t \right\}$ such that $X$ is smooth along the general curve $D_t$ and such that 
\[((p-1)L-K_X) \cdot D_t >0. \]
Then for every point $x \in X$ there exists a rational curve $C_x$ passing through $x$ such that
\[L \cdot C_x \leq 2 \dim X \frac{L \cdot D_t}{((p-1)L-K_X) \cdot D_t}.\]
\end{thm}
As an application we deduce an effective vanishing for the $H^1$ of a positive line bundle on a klt del Pezzo surface:
\begin{thm}\label{vanishing}
Let $k$ be an algebraically closed field of characteristic $p>0$.
Let $X$ be a klt del Pezzo surface over $k$ and let $A$ be a big and nef Cartier divisor. Then
\begin{enumerate}
\item $H^1(X,\mathcal{O}_X(-A))=0$; 
\item If $p \geq 5$, then $H^1(X,\mathcal{O}_X(A))= 0$;
\item If $p=3$, then $H^1(X,\mathcal{O}_X(2A))= 0$;
\item If $p=2$, then $H^1(X, \mathcal{O}_X(4A))=0$. 
\end{enumerate}
\end{thm}
\begin{proof}
To prove (1), it is enough to show that $H^0(X, \mathcal{O}_X(A)) \neq 0$ by \cite[Proposition 3.3]{CT}. So denoting by $f \colon Y \rightarrow X$ the minimal resolution, we have
\[H^0(X, \mathcal{O}_X(A)) =H^0(Y,\mathcal{O}_Y(f^*A)). \]
Since $Y$ is a rational surface by Lemma \ref{l-rational-dP} we have $h^2(Y,\mathcal{O}_Y(f^*A))=h^0(Y,\mathcal{O}_Y(K_Y-f^*A))=0$ and therefore
\[h^0(Y, \mathcal{O}_Y(f^*A)) \geq 1+ \frac{1}{2}f^*A \cdot (f^*A-K_Y)=1+ \frac{1}{2} (A^2-K_X \cdot A) > 0. \] 
To prove (2), let us note that if $H^1(X,\mathcal{O}_X(A)) \neq 0$ we have $H^1(X,\mathcal{O}_X(K_X-A)) \neq 0$ by Serre duality. Let us define a Weil $\Q$-Cartier ample divisor
\[L:= A-K_X. \]
Considering a covering family $\left\{D_t \right\}$ of curves for $X$ belonging to a very ample linear system we have that
\[(p-1)(L\cdot D_t)-K_X \cdot D_t>0. \]
Therefore we can apply Theorem \ref{kollarbend}  for every point $x \in X$ we can find a curve $C_x$ passing through $x$ such that
\[L \cdot C_x \leq 4 \frac{L \cdot D_t}{(p-1)L\cdot D_t-K_X \cdot D_t} < \frac{4}{p-1}. \]
Moreover, if $x \in X$ is chosen to be generic we have that $A \cdot C_x \geq 1$ since $A$ is big Cartier divisor and therefore 
\[L \cdot C_x = A\cdot C_x -K_X \cdot C_x > 1. \]
Thus concluding that $p < 5$. \\
\indent In the case where $p=3$, we apply the same proof to $L=2A-K_X$ with the same notation to the curves $D_t$. In this case by Theorem \ref{kollarbend} we can find that for any point $x$ there exists a rational curve $C_x$ passing through $x$ such that
\[L \cdot C_x < \frac{4}{p-1}=2. \]
However choosing $x$ generic enough we have
\[ L \cdot C_x=2A \cdot C_x -K_X \cdot C_x >2, \]
thus getting a contradiction. The proof for the case $p=2$ is analogous.
\end{proof}

We conclude by discussing the special case where the linear system induced by $A$ is birational.
 
\begin{proposition}
Let $k$ be an algebraically closed field of characteristic $p>0$.
Let $(X, \Delta)$ be a log del Pezzo pair over $k$. 
Let $A$ be a big and nef Cartier divisor such that the linear system $|A|$ is base point free and birational onto the image. Then
\[H^1(X,\mathcal{O}_X(A))=0. \]
\end{proposition}

\begin{proof}
Let $f \colon Y \rightarrow X$ be the minimal resolution. We have for a certain effective boundary divisor $\Delta_Y$:
\[K_Y + \Delta_Y = f^* (K_X+\Delta). \]
Since klt surface singularities are rational, we have $R^1f_* \mathcal{O}_Y (f^*A) = R^1f_* \mathcal{O}_Y \otimes \mathcal{O}_X(A)=0$. Thus we deduce
\[H^i(Y, \mathcal{O}_Y( f^*A)) = H^i(X, \mathcal{O}_X(A)). \]
By hypothesis, there exists an integral curve $C \in |f^*A|$. The curve $C$ is Gorenstein with dualizing sheaf:
\[ \omega_C=\mathcal{O}_Y(K_Y+C) \otimes \mathcal{O}_C. \] 
Consider the following short exact sequence:
\[0 \rightarrow \mathcal{O}_Y \rightarrow \mathcal{O}_Y(C) \rightarrow \mathcal{O}_C(C) \rightarrow 0. \]
Since $Y$ is a smooth rational surface by Lemma \ref{l-rational-dP} we have $H^i(Y, \mathcal{O}_Y)=0$ for $i=1,2$. Thus taking the long exact sequence in cohomology we have
\[H^1(Y, \mathcal{O}_Y(f^*A)) = H^1(Y, \mathcal{O}_Y(C)) \simeq H^1(C, \mathcal{O}_C(C)).\]
Now, using Serre duality on $C$ we have
\[ H^1(C, \mathcal{O}_C(C)) \simeq H^1(C, \omega_C \otimes \mathcal{{O}}_C(-K_Y|_C)) ) \simeq H^0(C, \mathcal{O}_C(K_Y|_C) )^*. \]
It is easy to see that $K_Y|_C$ is an anti-ample divisor because
\[K_Y \cdot C=(f^*(K_X+\Delta)-\Delta_Y)\cdot f^*A=(K_X+\Delta) \cdot A - \Delta_Y \cdot f^*A <0. \]
Therefore
\[H^0(C, \mathcal{O}_C(K_Y|_C))=0, \]
thus concluding the proof.
\end{proof}

\end{document}